\documentclass[12pt]{article}
\usepackage{latexsym,amsfonts,amssymb}
\usepackage{indentfirst}

\def\y{\begin{eqnarray*}}
\def\bd{\begin{description}}
\def\ey{\end{eqnarray*}}
\def\ebd{\end{description}}

\def\R{\mathbb{R}}

\begin{document}

\setlength{\baselineskip}{18pt}

\title{ Positive periodic solutions for abstract evolution equations with delay
\thanks{Research supported by NNSF of China (11261053) and NSF of Gansu Province (1208RJZA129).
}}
\vskip3mm
\author{Qiang Li $^{1}$\thanks{
 \emph{Corresponding author.}E-mail address: lznwnuliqiang@126.com (Q. Li), liyx@nwnu.edu.cn (Y. Li).,weimei@sxnu.edu.cn(M. Wei)}, Yongxiang Li $^{2}$, Mei Wei $^{1}$
}

\date{\small
\begin{flushleft}
$1.$ Department of Mathematics, Shanxi Normal University, Linfen 041000, Peoples's Republic of China,
\vskip1mm
\noindent
$2.$ Department of Mathematics, Northwest Normal University, Lanzhou 730070, Peoples's Republic of China
\end{flushleft}
}

\maketitle

\begin{abstract}
\setlength{\baselineskip}{16pt}

In this paper, we discuss the existence and asymptotic stability of
the positive periodic mild solutions for the abstract evolution
equation with delay in an ordered Banach space $E$,
  $$u'(t)+Au(t)=F(t,u(t),u(t-\tau)),\ \ \ \ t\in\R,$$
  where $A:D(A)\subset E\rightarrow E$ is a closed linear operator
   and $-A$ generates a positive $C_{0}$-semigroup $T(t)(t\geq0)$,
    $F:\R\times E\times E\rightarrow E$ is
     a continuous mapping which is $\omega$-periodic in $t$.
  Under the ordered conditions on the nonlinearity $F$ concerning the growth exponent
   of the semigroup $T(t)(t\geq0)$ or the first eigenvalue of the operator $A$,
   we obtain the existence and asymptotic  stability  of
   the positive $\omega$-periodic mild solutions by applying
    operator semigroup theory. In the end, an example
     is given to illustrate the applicability of our abstract
results.

\vspace{8pt}

\vskip3mm
\noindent {\bf Key Words:\ } Evolution equations with delay; Positive periodic solutions; Existence and uniqueness; Asymptotic stability; Positive $C_{0}$-semigroup
\vskip3mm
\noindent {\bf  MR(2010) Subject Classification:\ } 34K30; 47H07; 47H08
\vskip5mm
\end{abstract}

\section{Introduction and main results}

The theory of partial differential equation with delay has
 extensive physical background and realistic mathematical model,
  and it has undergone a rapid development in the last fifty years.
   Such equations are often more realistic to describe natural
phenomena than those without delay(see \cite{Hale1993,Wu1996}).

The problems concerning periodic solutions of partial
 differential equations with delay are an important area
  of investigation since they can take into account
   seasonal fluctuations occurring in the phenomena
   appearing in the models, and have been studied by
    some researchers in recent years.
The existence and asymptotic stability
 of periodic solutions of evolution
equation with delay have attracted much attention,
 see \cite{Burton1991,Xiang1992,Liu1998,Liu2000,Liu2003, Li2011,Kpoumie,Wang2016,Liang2017A,Liang2017D}.

Specilly, in \cite{Li2011}, by using analytic semigroups theory
and an integral inequality with
delays, Li discussed the time periodic solution for the evolution
equation with multiple delays in a Hilbert space $H$
$$u'(t)+Au(t)=F(t,u(t),u(t-\tau_{1}),\cdots,u(t-\tau_{n}))),\ \ \ \ t\in\R,\eqno(1.1)$$
where $A : D(A)\subset H\rightarrow H$ is a positive definite selfadjoint operator,
having compact resolvent and  the first eigenvalue  $\lambda_{1}> 0$,
$F : \R \times H^{n+1} \rightarrow H$ is a
nonlinear mapping which is $\omega$-periodic in $t$,
and $\tau_{1},\tau_{2},\cdots,\tau_{n}$ are positive constants which
denote the time delays. Under the following assumptions
\vskip1mm
\indent
(F1)
$\|F(t,v_0, v_1,\cdots , v_n)\|\leq \sum_{i=0}^{n} \beta_i\|v_i\|+K, t\in \R, (v_{0},\cdots v_n) \in H^{n+1}$,
\vskip1mm
\indent (F2) $\sum_{i=0}^{n} \beta_i<\lambda_{1}$,
\vskip1mm
\indent  (F3) $ \|F(t,v_0, v_1, \cdots , v_n)-F(t,w_0,w_1,\cdots, w_n)\|\leq \sum_{i=0}^{n} \beta_i\|v_i-w_{i}\|$,
\vskip1mm
\noindent  the author obtained the existence and uniqueness of
time $\omega$-periodic solutions to Eq. (1.1),
where  $\beta_{0},\beta_{1},\cdots,\beta_{n}$ and $K$ are positive constants. Moreover, strengthening the condition (F2) as follow
\vskip1mm
\indent (F2$^{*}$) $\beta_{0}+\sum_{i=1}^{n}e^{\lambda_{1}\tau_{i}\beta_{i}}<\lambda_{1}$,
\vskip1mm
\noindent
the unique time periodic solution was asymptotically stable.
However, because of the limitation of the research space and
the particularity of the operator, the results of the research
are not universal, and sometimes the conditions (F1) and (F3)
are not easy to verify in applications.

In \cite{Kpoumie}, Kpoumi\`{e} et al discussed the existence of periodic
solutions for the following nonautonomous
partial functional differential equation with  delay
$$u'(t)=A(t)u(t)+L(t,u_{t})+F(t,u_{t}),\ \ \ \ t\geq0 \eqno(1.2)$$
in a Banach space $X$, where $(A(t))_{t\geq0}$ is a family of linear operators on $X$,
 $L$ and $F$ are given continuous mappings and $\omega$-periodic
with respect to the first argument, the history $u_t$,
for $t\geq 0$, is defined from $(-\infty, 0]$ to $X$ by
$$u_{t}(s)=u(t+s),\ \ \ \ \  s\in (-\infty,0]. $$
 By using Massera's approach and  fixed point for multivalued maps,
 they proved the existence of an $\omega$-periodic solution.

 Recently, In \cite{Wang2016}, under suitable assumptions,
 such as the ultimate boundedness of the solutions of equations,
 Wang and Zhu established a theorem on periodic solutions to equations of this kind by using the Horn fixed-point theorem. In \cite{Liang2017A,Liang2017D},
 Liang et al also studied  nonautonomous  evolutionary
 equations with time delay and impulsive. Under the nonlinear term
 satisfying continuous and Lipschitzian, the proved the existence theorem for periodic mild solutions to the nonautonomous delay evolution equations by Horn's
 fixed point theorem or Sadovskii's Fixed Point Theorem.
However, in all these works, the key assumption or process of prior
boundedness of solutions was employed.

In many practice models, such as heat conduction equation,
neutron transport equation, reaction diffusion equation, etc.,
only positive periodic solutions are significant. In \cite{Li2005x},
the existence and uniqueness of positive periodic mild solutions for the evolution
equation without delay
$$u'(t)+Au(t)=F(t,u(t)),\ \ \ \ t\in\R,\eqno(1.3)$$
are obtained in an ordered Banach space $E$, where $-A$ is
 the infinitesimal generator of a positive $C_{0}$-semigroup,
 $F:\R\times E\to E$ is a continuous mapping which is $\omega$-periodic in $t$.
Recently, under the ordered
conditions on the nonlinearity $F$,
the existence and asymptotic stability of  positive
$\omega$-periodic mild solutions for the evolution
equation (1.3) have been obtained by applying operator semigroup theory,
monotone iterative technique and some fixed point theorems in an ordered
Banach space $E$,  see \cite{QLi2015}.
However, to the best of our knowledge, there are few papers to study
the existence and asymptotic stability of  positive
$\omega$-periodic solutions for the evolution
equation with delay. Furthermore, for the abstract evolution equation
without delay, the periodic solutions have been discussed by
more authors, see \cite{Amann76,Amann78,Li1996,Li98,Ahmed01, Li2009,ShiLiSun11,Yang}
and references therein.

Motivated by the papers mentioned above,  by means of
operator semigroup theory and some fixed point theorems,
we will use a completely different method to improve and
extend the results mentioned above,
which will make up the research in this area blank.

Our discussion will be made in the framework of ordered Banach spaces.
Let $E$ be an ordered Banach space $E$, whose positive cone $K$ is
normal cone with normal constant $N$. Let $A:D(A)\subset E\rightarrow E$
is a closed linear operator and $-A$ generates a positive $C_{0}$-semigroup
$T(t)(t\geq0)$ in $E$, the nonlinear function $F:\R\times E\times E\rightarrow E$
is a continuous mapping  and for every $x,y\in K$,
 $F(t,x,y)$ is $\omega$-periodic in $t$.
In this paper, we consider the following abstract evolution equation with delay
$$u'(t)+Au(t)=f(t,u(t),u(t-\tau)),\qquad t\in\R, \eqno(1.4)$$
We will study the existence and asymptotic stability of positive
$\omega$-periodic mild solutions for (1.4) under
some new conditions by applying the Leray-Schauder
fixed point theorem in an ordered
Banach space $E$. More precisely, the nonlinear term satisfies
order conditions concerning the growth exponent of the
semigroup $T(t)(t\geq0)$ or the first eigenvalue of
the operator $A$.

For $C_{0}$-semigroup $T(t)(t\geq0)$, there exist $M>0$ and
$\gamma\in \R$ such that (see \cite{Pazy})
$$\|T(t)\|\leq Me^{\gamma t},\quad t\geq0.\eqno(1.5)$$
Let
$$\nu_{0}=\inf\{\gamma\in \R |\ \mathrm{There\ exists}\ M>0\ \mathrm{ such \ that}\  \|T(t)\|\leq Me^{\gamma t},\ \forall t\geq0\},$$
 then $\nu_{0}$ is called the growth exponent of the semigroup $T(t)(t\geq0)$.
 Furthermore, $\nu_{0}$ can be also obtained by the following formula
$$\nu_{0}=\limsup\limits_{t\rightarrow +\infty}\frac{\ln\|T(t)\|}{t}.$$
If $C_{0}$-semigroup $T(t)$ is continuous in the
uniform operator topology for every $t>0$ in $E$,
it is well known that $\nu_{0}$ can also be determined by $\sigma(A)$ (see \cite{Triggiani})
$$\nu_{0}=-\inf\{\mathrm{Re} \lambda |\ \lambda\in\sigma(A)\},\eqno(1.6)$$
where $-A$ is the infinitesimal generator of
$C_{0}$-semigroup $T(t)(t\geq0)$. We know that
$T(t)(t\geq0)$ is continuous in the uniform
operator topology for $t>0$ if $T(t)(t\geq0)$ is compact semigroup (see \cite{Triggiani}).

For the abstract delay evolution equation (1.4),
we obtain the following results:
\vskip3mm
\noindent\textbf{Theorem 1.1.} \emph{Let $-A$ generate
an exponentially stable positive compact semigroup $T(t)(t\geq0)$ in $E$,
that is $\nu_{0}<0$. Assume that $F:\R\times K\times K\rightarrow K$
is a continuous mapping  which is $\omega$-periodic in $t$.
If the following condition
\vskip2mm
\noindent(H1) there are positive constants $C_{1},C_{2}$
satisfying $C_{1}+C_{2}\in(0,|\nu_{0}|)$ and a function $h_{0}\in C_{\omega}(\R,K)$ such that
$$F(t,x,y)\leq C_{1}x+C_{2}y+h_{0}(t),\quad t\in\R,\ x,y\in K,$$
\vskip2mm
\noindent holds, then Eq.$(1.4)$ has at least one
positive $\omega$-periodic mild solution $u$.
}
\vskip3mm
\noindent\textbf{Theorem 1.2.} \emph{ Let $-A$ generate an
exponentially stable positive compact semigroup $T(t)(t\geq0)$ in $E$.
Assume that $F:\R\times K\times K\rightarrow K$ is
a continuous mapping  which is $\omega$-periodic in $t$.
If the following condition
\vskip2mm
\noindent(H2)  there are positive constants $C_{1},C_{2}$
satisfying $C_{1}+C_{2}\in(0,|\nu_{0}|)$, such that for any $x_{i},y_{i}\in K(i=1,2)$ with $x_{1}\leq x_{2}, y_{1}\leq y_{2}$,
 $$F(t,x_{2},y_{2})-F(t,x_{1},y_{1})\leq C_{1}(x_{2}-x_{1})+C_{2}(y_{2}-y_{1}),\ \qquad t\in\R, $$
\noindent holds, then Eq. $(1.4)$ has a unique
positive $\omega$-periodic mild solution $u$.
}

Now, we strengthen the condition (H2) in Theorem 1.2,
we can obtain the following
asymptotic stability result of the periodic solution:

\vskip3mm
\noindent\textbf{Theorem 1.3.} \emph{Let $-A$ generate an exponentially
stable positive compact semigroup $T(t)(t\geq0)$ in $E$.
Assume that $F:\R\times K\times K\rightarrow K$ is a continuous mapping
which is $\omega$-periodic in $t$. If the following condition
\vskip2mm
\noindent(H3)  there are positive constants $C_{1},C_{2}$ satisfying $C_{1}+C_{2}e^{-\nu_{0}\tau}\in(0,|\nu_{0}|)$, such that for any $x_{i},y_{i}\in K(i=1,2)$ with $x_{1}\leq x_{2}, y_{1}\leq y_{2}$,
$$F(t,x_{2},y_{2})-F(t,x_{1},y_{1})\leq C_{1}(x_{2}-x_{1})+C_{2}(y_{2}-y_{1}),\ \qquad t\in\R, $$
\noindent holds, then the unique positive $\omega$-periodic
mild solution of Eq.(1.4) is globally asymptotically stable.
}

Furthermore, we assume that the  positive cone $K$ is regeneration cone.
By the characteristic of positive semigroups (see \cite{Li1996}), for sufficiently large
$\lambda_{0}>-\inf\{Re \lambda| \lambda\in\sigma(A)\}$,
we have that $\lambda_{0}I+A$ has positive bounded
inverse operator $(\lambda_{0}I+A)^{-1}$.
Since $\sigma(A)\neq \emptyset $, the spectral radius $r((\lambda_{0}I+A)^{-1})=\frac{1}{dist(-\lambda_{0},\sigma(A))}>0$.
By the famous Krein-Rutmann theorem,
$A$ has the first eigenvalue $\lambda_{1}$,
which has a positive eigenfunction $e_{1}$, and
$$\lambda_{1}=\inf \{Re\lambda |\ \lambda\in \sigma(A)\},\eqno(1.7)$$
that is $\nu_{0}=-\lambda_{1}$. Hence, by Theorem 1.1,
Theorem 1.2 and Theorem 1.3, we have the following results.

\vskip3mm
\noindent\textbf{Corollary 1.4}\emph{ Let  $-A$ generate an exponentially
stable positive compact semigroup $T(t)(t\geq0)$ in $E$.
Assume that $F:\R\times K\times K\rightarrow K$ is
a continuous mapping  which is $\omega$-periodic in $t$.
If the following condition
\vskip2mm
\noindent(H1$'$) there are positive constants $C_{1},C_{2}$
satisfying $C_{1}+C_{2}\in(0,\lambda_{1})$ and a
function $h_{0}\in C_{\omega}(\R,K)$ such that
$$F(t,x,y)\leq C_{1}x+C_{2}y+h_{0}(t),\quad t\in\R,\ x,y\in K,$$
\vskip2mm
\noindent holds, then Eq.$(1.4)$ has at least one
positive $\omega$-periodic mild solution $u$.
}
\vskip3mm
\noindent\textbf{Corollary 1.5} \emph{ Let $-A$ generate an exponentially
stable positive compact semigroup $T(t)(t\geq0)$ in $E$.
Assume that $F:\R\times K\times K\rightarrow K$ is a
continuous mapping  which is $\omega$-periodic in $t$.
If the following condition
\vskip2mm
\noindent(H2$'$)  there are positive constants $C_{1},C_{2}$
satisfying $C_{1}+C_{2}\in(0,\lambda_{1})$, such that for any $x_{i},y_{i}\in K(i=1,2)$ with $x_{1}\leq x_{2}, y_{1}\leq y_{2}$,
 $$F(t,x_{2},y_{2})-F(t,x_{1},y_{1})\leq C_{1}(x_{2}-x_{1})+C_{2}(y_{2}-y_{1}),\ \qquad t\in\R, $$
\noindent holds, then Eq. $(1.4)$ has a
unique positive $\omega$-periodic mild solution $u$.
}

\vskip3mm
\noindent\textbf{Corollary 1.6} \emph{ Let $-A$ generate an exponentially
stable positive compact semigroup $T(t)(t\geq0)$ in $E$.
Assume that $F:\R\times K\times K\rightarrow K$
is a continuous mapping  which is $\omega$-periodic in $t$.
If the following condition
\vskip2mm
\noindent(H3$'$)  there are positive constants $C_{1},C_{2}$ satisfying $C_{1}+C_{2}e^{\lambda_{1}\tau}\in(0,\lambda_{1})$, such that for any $x_{i},y_{i}\in K(i=1,2)$ with $x_{1}\leq x_{2}, y_{1}\leq y_{2}$,
$$F(t,x_{2},y_{2})-F(t,x_{1},y_{1})\leq C_{1}(x_{2}-x_{1})+C_{2}(y_{2}-y_{1}),\ \qquad t\in\R, $$
\noindent holds, then the unique positive $\omega$-periodic
 mild solution of Eq. (1.4) is globally asymptotically stable.
}
\vskip3mm
\noindent\textbf{Remark 1.7.}\emph{ In Corollary 1.4 and Corollary 1.5, since $\lambda_{1}$
is the first eigenvalue of $A$, the
condition $C_{1}+C_{2}<\lambda_{1}$ in (H1$'$) and (H2$'$)
cannot be extended to $C_{1}+C_{2}\leq\lambda_{1}$.
Otherwise,  periodic problem (1.4) does not
always have a mild solution.
For example, $F(t,x,y)=\frac{\lambda_{1}}{2}x+\frac{\lambda_{1}}{2}y$.}
\vskip3mm
\noindent\textbf{Remark 1.8.}\emph{ It is clear that our results
can also be extended to  the evolution equation with multiple delays (1.1).
In this case, the conditions (H1$'$) and (H3$'$) have improved the conditions
 (F1) and (F3), and our conditions are  easy to verify in applications.
 Hence, our results of the positive periodic solutions,
  improve and generalize the results in \cite{Li2011}.
 On the other hand, we delete the Lipschitz conditions on nonlinearity.
 In this case, the prior estimate of solutions are not employed. Thus, compared with the existence results in \cite{Kpoumie,Liang2017A,Liang2017D}, our conclusions are new in some respects in some respects.
 }

\vskip3mm

The paper is organized as follows. Section 2 provides the
definitions and preliminary results to be used in theorems
stated and proved in the paper. The proofs of Theorems 1.1-1.3
are based on positive $C_{0}$-semigroups theory, Leray-Schauder
fixed point theorem and an integral inequality of Bellman, which
will be given in Section 3. In the last section, we give an example
 to illustrate the applicability of the abstract results.

\section{Preliminaries}

In this section, we introduce some notions, definitions,
and preliminary facts which are used through this paper.

Let $J$ denote the infinite interval $[0,+\infty)$ and $h:J\rightarrow E$,
consider the initial value problem of the linear evolution equation
$$\left\{\begin{array}{ll}
u'(t)+Au(t)=h(t),\ t\in J,\\[8pt]
u(0)=x_{0}
 \end{array} \right.\eqno(2.1)$$
 It is well known \cite[Chapter 4, Theorem 2,9]{Pazy},
 when $x_{0}\in D(A)$ and $h\in C^{1}(J,E)$,
 the initial value problem (2.1) has a unique classical solution
 $u\in C^{1}(J,E)\cap C(J,E_{1})$ expressed by
 $$u(t)=T(t)x_{0}+\int^{t}_{0}T(t-s)h(s)ds, \eqno(2.2)$$
where $E_{1}=D(A)$ is
Banach space  with the graph norm $\|\cdot\|_{1} = \|\cdot\|+\|A\cdot\|$.
Generally, for $x_{0}\in E$ and $h\in C(J,E)$,
the function $u$ given by (2.2) belongs to $C(J,E)$ and it is
called a mild solution of the linear evolution equation (2.1).

Let $C_{\omega}(\R,E)$ denote the Banach space $ \{u\in C(\R,E)|~u(t+\omega)=u(t),\ t\in \R\} $
endowed the maximum norm $\|u\|_{C}=\max_{t\in[0,\omega]}\|u(t)\|$.
Evidently, $C_{\omega}(\R,E)$ is also an order Banach space with the partial order $``\leq"$
induced by the positive cone $K_{C}=\{u\in C_{\omega}(\R,E)|\  \ u(t)\geq \theta,\ t\in \R\}$
and $K_{C}$ is also normal with the normal constant $N$.

Given $h\in C_{\omega}(\R,E)$, for the following linear evolution
 equation corresponding to Eq.(1.4)
$$u'(t)+Au(t)=h(t),\quad t\in \R,\eqno(2.3)$$
we have the following result.
\vskip3mm
\noindent\textbf{Lemma 2.1.}(\cite{Li2005x}) \emph{ If $-A$
generates an exponentially stable positive
$C_{0}$-semigroup $T(t)(t\geq0)$ in $E$, then for $h\in C_{\omega}(\R,E)$,
the linear evolution equation (2.3) exists a unique positive $\omega$-periodic mild
solution $u$, which can be expressed by
  $$u(t)=(I-T(\omega))^{-1}\int_{t-\omega}^{t}T(t-s)h(s)ds:=(Ph)(t),\eqno(2.4)$$
   and the solution operator $P:C_{\omega}(\R,E)\rightarrow C_{\omega}(\R,E)$
    is a positive bounded linear operator
     with the spectral radius $r(P)\leq\frac{1}{|\nu_{0}|}$.}
\vskip3mm
\noindent\textbf{Proof.} For any $\nu\in (0,|\nu_{0}|)$,
there exists $M>0$ such that
$$\|T(t)\|\leq Me^{-\nu t}\leq M, \quad t\geq0.\eqno(2.5)$$
In $E$, define the equivalent norm $|\cdot|$ by
$$|x|=\sup\limits_{t\geq0}\|e^{\nu t}T(t)x\|,$$
then $\|x\|\leq |x|\leq M\|x\|$.
By $|T(t)|$ we denote the norm of $T(t)$
in $(E,|\cdot|$), then for $t\geq0$, it is easy to obtain that
$|T(t)|<e^{-\nu t}$.
Hence, $(I-T(\omega))$ has bounded inverse operator
$$(I-T(\omega))^{-1}=\sum_{n=0}^{\infty}T(n\omega),\eqno(2.6)$$
and its norm satisfies
$$|(I-T(\omega))^{-1}|\leq \frac{1}{1-|T(\omega)|}\leq \frac{1}{1-e^{-\nu \omega}}.\eqno(2.7)$$
Set $$x_{0}=(I-T(\omega))^{-1}\int^{\omega}_{0}T(t-s)h(s)ds:=Bh,\eqno(2.8)$$
 then the mild solution $u(t)$  of
 the linear initial value problem  (2.1)
given by (2.2) satisfies the periodic boundary
 condition $u(0)=u(\omega)=x_{0}$.
For $t\in \mathbb{R}^{+}$, by (2.2) and the
 properties of the semigroup $T(t)(t\geq0)$,
we have
\begin{eqnarray*}
u(t+\omega)&=&T(t+\omega)u(0)+\int_{0}^{t+\omega}T(t+\omega-s)h(s)ds\\[8pt]
&=&T(t)\Big(T(\omega)u(0)+\int_{0}^{\omega}T(\omega-s)h(s)ds\Big)+\int^{t}_{0}T(t-s)h(s-\omega)ds\\[8pt]
&=&T(t)u(0)+\int^{t}_{0}T(t-s)h(s)ds\ =u(t).
\end{eqnarray*}
Therefore, the $\omega$-periodic extension of $u$
on $\R$ is a unique $\omega$-periodic
mild solution of Eq.(2.3). By (2.2) and (2.8),
the $\omega$-periodic mild solution
can be  expressed by
\begin{eqnarray*}\qquad\qquad\quad
u(t)&=&T(t)B(h)+\int^{t}_{0}T(t-s)h(s)ds\\[8pt]
&=&(I-T(\omega))^{-1}\int^{t}_{t-\omega}T(t-s)h(s)ds:=(Ph)(t).\quad \qquad\qquad\ (2.9)
\end{eqnarray*}
Evidently, by the positivity
of semigroup $T(t)(t\geq0)$, we can obtain that
 $P:C_{\omega}(\R,E)\rightarrow C_{\omega}(\R,E)$
is a positive bounded linear operator. By (2.7) and (2.9), we have
\begin{eqnarray*}
&&|(Ph)(t)|\leq|(I-T(\omega))^{-1}|\int_{t-\omega}^{t}|T(t-s)h(s)|ds\\[8pt]
&\leq&\frac{1}{1-e^{-\nu \omega}}\int_{t-\omega}^{t}e^{-\nu(t-s)}|h|_{C}ds\\[8pt]
&\leq&\frac{1}{\nu}|h|_{C},
\end{eqnarray*}
which implies that $|P|\leq \frac{1}{\nu}$.
Therefore, $r(P)\leq |P|\leq\frac{1}{\nu}$.
Hence, by the arbitrary of $\nu\in(0,|\nu_{0}|)$,
we have $r(P)\leq\frac{1}{|\nu_{0}|}$.
This completes the proof of Lemma 2.1. $\Box$

In the proof of our main results, we also need the following results.

\vskip3mm
\noindent\textbf{Lemma 2.2.}( Leray-Schauder fixed point
theorem \cite{Deimling})\emph{ Let $\Omega$  be convex
 subset of Banach space $E$ with $\theta\in \Omega$, and
 let $Q: \Omega\rightarrow \Omega$ be compact  operator.
 If the set $\{u\in \Omega |\ u=\eta Qu,\ 0<\eta<1\}$ is bounded,
 then $Q$ has a fixed point in $\Omega$.
}

\section{Proof of the main results}

\vskip3mm
\noindent\textbf{Proof of  Theorem 1.1.} Evidently, the normal
cone $K_{C}$  is a convex subset of Banach space $C_{\omega}(\R,E)$ and $\theta\in K_{C}$.
Consider the operator $Q$ defined by
$$Qu=(P\circ \mathcal{F})(u),\eqno(3.1)$$
where
$$\mathcal{F}(u)(t):=F(t,u(t),u(t-\tau)),\quad u\in K_{C}.\eqno(3.2)$$
From the positivity of semigroup of $T(t)(t\geq0)$
and the conditions of Theorem 1.1, it is easy to see
that $Q:K_{C}\rightarrow K_{C} $ is well defined.
From (3.1) and (3.2), it follows that
$$(Qu)(t)=(I-T(\omega))^{-1}\int_{t-\omega}^{t}T(t-s)F(s,u(s),u(s-\tau))ds,\quad t\in \R.\eqno(3.3)$$
By the definition $P$, the positive $\omega$-periodic mild solution
of Eq.(1.4) is equivalent to
the fixed point of the operator $Q$. In the following,
we will prove $Q$ has a fixed point by
applying the famous Leray-Schauder fixed point theorem.

At first, we prove that $Q$ is continuous on $K_{C}$.
Let $\{u_{m}\}\subset K_{C}$ be a sequence
such that $u_{m}\rightarrow u\in K_{C} $ as $m\rightarrow \infty$, so for every $t\in \R$, $\lim\limits_{m\rightarrow\infty}u_{m}(t)=u(t)$.
Since $F:\R\times K^{n+1}\rightarrow K$ is continuous,
then for  every $t\in \R$, we get
$$F(t,u_{m}(t),u_{m}(t-\tau))\rightarrow F(t,u(t),u(t-\tau)),\quad m\rightarrow \infty.\eqno(3.4) $$
By (3.3) and the Lebesgue dominated convergence theorem,
for every $t\in \R$, we have
\begin{eqnarray*}\qquad
& &\|(Qu_{m})(t)-(Qu)(t)\|\\[8pt]
&=&\Big\|(I-T(\omega))^{-1}\Big(\int_{t-\omega}^{t}T(t-s)F(s,u_{m}(s),u_{m}(s-\tau))ds\\[8pt]
&&-
\int_{t-\omega}^{t}T(t-s)F(s,u(s),u(s-\tau))ds\Big)\Big\|\\[8pt]
&\leq&\|(I-T(\omega))^{-1}\|\cdot\int_{t-\omega}^{t}\|T(t-s)\|\cdot\|F(s,u_{m}(s),u_{m}(s-\tau))\\[8pt]
&&- F(s,u(s),u(s-\tau))\|ds\\[8pt]
&\leq&CM\cdot\int_{t-\omega}^{t}\|F(s,u_{m}(s),u_{m}(s-\tau))- F(s,u(s),u(s-\tau))\|ds, \ \qquad(3.5)
\end{eqnarray*}
where $C=\|(I-T(\omega))^{-1}\|$.
Therefore, we can conclude that
$$\|Qu_{m}-Qu\|\rightarrow 0,\quad m\rightarrow \infty.\eqno(3.6)$$
Thus, $Q:K_{C}\rightarrow K_{C}$ is continuous.

Subsequently, we show that $Q$ maps every bounded
set in $K_{C}$ into a bounded
set. For any $R>0$, let
$$\overline{\Omega}_{R}:=\{u\in K_{C}| \ \|u\|_{C}\leq R\}. \eqno(3.7)$$
For each $u\in \overline{\Omega}_{R}$, from the continuity of $F$,
we know that there exists $M_{1}>0$ such that
$$\|F(t,u(t),u(t-\tau))\|\leq M_{1}, \quad t\in \R,\eqno(3.8)$$
hence, we get
\begin{eqnarray*}
\|(Qu)(t)\|&=&\|(I-T(\omega))^{-1}\int_{t-\omega}^{t}T(t-s)F(s,u(s),u(s-\tau))ds\|\\[8pt]
&\leq&\|(I-T(\omega))^{-1}\|\int_{t-\omega}^{t}\|T(t-s)\|\cdot \|F(s,u(s),u(s-\tau))\|ds\\[8pt]
&\leq& CM \int_{t-\omega}^{t}M_{1}ds\\[8pt]
&\leq& CMM_{1}\omega :=\overline{R}.
\end{eqnarray*}
Therefore, $Q(\overline{\Omega}_{R})$ is bounded.

Next, we  demonstrate that $ Q(\overline{\Omega}_{R})$ is equicontinuous.
For every $u\in\overline{\Omega}_{R}$, by the periodicity of $u$,
we only consider it on $[0,\omega]$.
Set $ 0\leq t_{1}<t_{2}\leq\omega$, we get that
\begin{eqnarray*}
&&Qu(t_{2})-Qu(t_{1})\\[8pt]
&=&(I-T(\omega))^{-1}\int^{t_{2}}_{t_{2}-\omega}T(t_{2}-s)F(s,u(s),u(s-\tau))ds\\[8pt]
&~&-(I-T(\omega))^{-1}\int^{t_{1}}_{t_{1}-\omega}T(t_{1}-s)F(s,u(s),u(s-\tau))ds\\[8pt]
&=&(I-T(\omega))^{-1}\int^{t_{1}}_{t_{2}-\omega}(T(t_{2}-s)-T(t_{1}-s))F(s,u(s),u(s-\tau))ds\\[8pt]
&~&-(I-T(\omega))^{-1}\int^{t_{2}-\omega}_{t_{1}-\omega}T(t_{1}-s)F(s,u(s),u(s-\tau))ds\\[8pt]
&~&+(I-T(\omega))^{-1}\int^{t_{2}}_{t_{1}}T(t_{2}-s)F(s,u(s),u(s-\tau))ds\\[8pt]
&:=&I_{1}+I_{2}+I_{3},
\end{eqnarray*}
It is clear that
$$\|Qu(t_{2})-Qu(t_{1})\|\leq \|I_{1}\|+\|I_{2}\|+\|I_{3}\|.\eqno(3.9)$$
Now, we only need to check $\|I_{i}\|$
tend to $0$  independently of $u\in\overline{\Omega}_{R}$
when $t_{2}-t_{1}\rightarrow 0,i=1,2,3$.
From the definition of $I_{i}$, we can easily see
\begin{eqnarray*}
\|I_{1}\|&\leq&C\cdot\int^{t_{1}}_{t_{2}-\omega}\|(T(t_{2}-s)-T(t_{1}-s))\|\cdot\|F(s,u(s),u(s-\tau))\|ds\\[8pt]
&\leq & CM_{1}\int^{t_{1}}_{t_{2}-\omega}\|(T(t_{2}-s)-T(t_{1}-s))\| ds\\[8pt]
&\rightarrow&0, \ \mathrm{as} \ t_{2}-t_{1}\rightarrow 0,\\[8pt]
\|I_{2}\|&\leq&C\cdot\int^{t_{2}-\omega}_{t_{1}-\omega}\|(T(t_{1}-s))\|\cdot\|F(s,u(s),u(s-\tau))\|ds\\[8pt]
&\leq &CMM_{1}(t_{2}-t_{1})\\[8pt]
&\rightarrow&0,  \ \mathrm{as} \ t_{2}-t_{1}\rightarrow 0,\\[8pt]
\|I_{3}\|&\leq&C\cdot\int^{t_{2}}_{t_{1}}\|(T(t_{2}-s))\|\cdot\|F(s,u(s),u(s-\tau))\|ds\\[8pt]
&\leq &CMM_{1}(t_{2}-t_{1})ds \\[8pt]
&\rightarrow&0,  \ \mathrm{as} \ t_{2}-t_{1}\rightarrow 0.
\end{eqnarray*}
As a result, $\| Qu(t_{2}) -Qu(t_{1})\|$ tends to $0$
 independently of $u\in \overline{\Omega}_{R}$
as $t_{2}- t_{1}\rightarrow0$,
which means that $Q(\overline{\Omega}_{R})$ is equicontinuous.

Now, we prove that $(Q\overline{\Omega}_{R})(t)$ is relatively
 compact in $K$ for all $t\in \R$.
We define a set $(Q_{\varepsilon}\overline{\Omega}_{R})(t)$ by
$$(Q_{\varepsilon}\overline{\Omega}_{R})(t):=\{(Q_{\varepsilon}u)(t)|\ u\in \overline{\Omega}_{R},\ 0<\varepsilon<\omega,\ t\in \R\}, \eqno(3.10)$$
where
\begin{eqnarray*}
(Q_{\varepsilon}u)(t)&=&(I-T(\omega))^{-1}\int_{t-\omega}^{t-\varepsilon}T(t-s)F(s,u(s),u(s-\tau))ds\\[8pt]
&=&T(\varepsilon)(I-T(\omega))^{-1}\int_{t-\omega}^{t-\varepsilon}T(t-s-\varepsilon)F(s,u(s),u(s-\tau))ds.
\end{eqnarray*}
Then the set $(Q_{\varepsilon}\overline{\Omega}_{R})(t)$ is
relatively compact in $K$ since the operator $T(\varepsilon)$ is compact in $K$.
For any $u\in \overline{\Omega}_{R}$ and $t\in \R$,
from the following inequality
\begin{eqnarray*}\qquad\qquad\qquad
&&\|Qu(t)-Q_{\varepsilon}u(t)\|\\[8pt]
&=&\Big\|(I-T(\omega))^{-1}\Big(\int_{t-\omega}^{t}T(t-s)F(s,u(s),u(s-\tau))ds\\[8pt]
&&-\int_{t-\omega}^{t-\varepsilon}T(t-s)F(s,u(s),u(s-\tau))ds\Big)\Big\|\\[8pt]
&\leq& C\int^{t}_{t-\varepsilon}\|T(t-s)F(s,u(s),u(s-\tau))\|ds\\[8pt]
&\leq& CMM_{1}\varepsilon,\ \ \quad \qquad\qquad\qquad\qquad\qquad\qquad\quad\qquad\qquad\qquad(3.11)\end{eqnarray*}
one can obtain that the set $(Q\overline{\Omega}_{R})(t)$
is relatively compact in $K$ for all $t\in \R$.

Thus, the Arzela-Ascoli theorem guarantees that
$Q:K_{C}\rightarrow K_{C}$ is a compact operator.

Finally, we prove the set
$\Lambda(Q):=\{u\in K_{C}|\ u=\eta Qu,\ \forall \ 0<\eta<1\}$ is bounded.
For every $u\in K_{C}$, by (3.2) and the condition (H1), we have
\begin{eqnarray*}\qquad\quad\qquad\quad \theta&\leq& \mathcal{F}(u)(t)=F(t,u(t),u(t-\tau))\\[8pt]
&\leq& C_{1}u(t)+C_{2}u(t-\tau)+h_{0}(t),\quad t\in \R.\qquad\quad \qquad\quad\qquad (3.12)\end{eqnarray*}
Define an operator $\mathcal{B}:C_{\omega}(\R,K)\to C_{\omega}(\R,K)$ as following:
$$\mathcal{B}u(t)=C_{1}u(t)+C_{2}u(t-\tau),\ \ \ t\in\R,\ u\in K_{C}. \eqno(3.13)$$
It is easy to see that $\mathcal{B}:C_{\omega}(\R,K)\to C_{\omega}(\R,K)$
 is a positive bounded linear operators satisfying
  $\|\mathcal{B}\|\leq C_{1}+C_{2}$.
  Let $u\in \Lambda(Q)$, then there is a constant $\eta\in(0,1)$
  such that $u=\eta Qu$. Therefore,
  by the definition of $Q$, Lemma 2.1 and (3.12), we have
\begin{eqnarray*}
\theta&\leq&u(t)=\eta (Qu)(t)<(Qu)(t)\\[8pt]
&=&P\circ \mathcal{F} (u)(t)\leq P(\mathcal{B}u(t)+h_{0}(t))\\[8pt]
&=&\mathcal{B}Pu(t)+Ph_{0}(t)<\mathcal{B}P\circ Q(t)+Ph_{0}(t)\\[8pt]
&\leq&cP(cPu(t)+Ph_{0}(t))+Ph_{0}(t)\\[8pt]
&=&\mathcal{B}^{2}P^{2}u(t)+\mathcal{B}P^{2}h_{0}(t)+Ph_{0}(t),
\end{eqnarray*}
inductively, we can see
$$u(t)\leq \mathcal{B}^{n}P^{n}u(t)+\mathcal{P}h_{0}(t),\quad n=1,2,\cdots,\eqno(3.14)$$
where, $\mathcal{P}=\mathcal{B}^{n-1}P^{n}+\mathcal{B}^{n-2}P^{n-1}+\cdots+\mathcal{B}P^{2}+P$
is a bounded linear operator, and there exists
a constant $M_{2}>0$ such that $\|\mathcal{P}\|\leq M_{2}$.
 Hence, by the normality of the cone $K_{C}$, we can see
\begin{eqnarray*}
\|u\|_{C}&<&N\|\mathcal{B}^{n}\|\cdot\|P^{n}\|\cdot\|u\|_{C}+M_{2}\|h_{0}\|_{C},\\[8pt]
&\leq&N(C_{1}+C_{2})^{n}\cdot\|P^{n}\|\cdot\|u\|_{C}+M_{2}\|h_{0}\|_{C}.\end{eqnarray*}
 From the spectral radius of Gelfand formula  $\lim\limits_{n\rightarrow\infty}\sqrt[n]{\|P^{n}\|}=r(P)=\frac{1}{|\nu_{0}|}$,
 and  the condition (H1), when $n$ is large enough, we get that $(C_{1}+C_{2})^{n}\cdot\|P^{n}\|<\frac{1}{N}$,
  then
$$\|u\|_{C}< \frac{M_{2}\|h_{0}\|_{C}}{1-N(C_{1}+C_{2})^{n}\cdot\|P^{n}\|},\eqno(3.15)$$
which implies that $\Lambda(Q)$ is bounded.
 By the Leray-Schauder fixed point theorem of compact operator,
 the operator $Q$ has at least one fixed point $u$ in $K_{C}$,
 which is a positive $\omega$-periodic mild solution of
 the delay evolution equation (1.4).
 This completes the proof of Theorem 1.1. \hfill$\Box$
\vskip3mm

\noindent\textbf{Proof of  Theorem 1.2.} From the condition (H3),
it is easy to see that the condition (H1) holds.
Hence by Theorem 1.1, Eq.(1.4)
 has positive $\omega$-periodic mild solutions.
 Let $u_{1},u_{2}\in K_{C}$ be the positive $\omega$-periodic solutions of Eq.(1.4),
 then they are the fixed points of the operator $Q=P\circ \mathcal{F}$.
 Let us assume $u_{1}\leq u_{2}$,
 by the definition of $\mathcal{F}$ and the condition (H3),
 for any $t\in\R$, we have
\begin{eqnarray*}
&&\mathcal{F}(u_{2})(t)-\mathcal{F}(u_{1})(t)\\[8pt]
&\leq&C_{1}(u_{2}(t)-u_{1}(t))
+C_{2}(u_{2}(t-\tau)-u_{1}(t-\tau))\\[8pt]
&=&\mathcal{B}(u_{2}(t)-u_{1}(t)),
\end{eqnarray*}
where $\mathcal{B}$ is defined by (3.13).
 Thus, we can obtain that
\begin{eqnarray*}
&&\theta\leq u_{2}(t)-u_{1}(t)=(Qu_{2})(t)-(Qu_{1})(t)\\[8pt]
&=&P((\mathcal{F}u_{2})(t)-(\mathcal{F}u_{1})(t))
\leq \mathcal{B}P(u_{2}(t)-u_{1}(t))\\[8pt]
&\leq&\cdots\leq \mathcal{B}^{n}P^{n}(u_{2}(t)-u_{1}(t)).
\end{eqnarray*}
By the normality of the cone $K_{C}$, we can see
$$\|u_{2}-u_{1}\|_{C}\leq N\|\mathcal{B}^{n}\|\cdot\|P^{n}\|\cdot\|u_{2}-u_{1}\|_{C},\eqno(3.16)$$
 From the proof of Theorem 1.1, when $n$ is large enough,
  $N\|\mathcal{B}^{n}\|\cdot\|P^{n}\|<1$, so $\|u_{2}-u_{1}\|_{C}=0$,
   it follows that $u_{2}\equiv u_{1}$.
   Thus, Eq.(1.4) has only one positive
    $\omega$-periodic mild solution.  \hfill$\Box$

\vskip3mm
In order prove Theorem 1.3, we need discuss
 the existence and uniqueness of
 the initial value problem of the nonlinear
  delay evolution equation (1.4).

Let $C([-\tau,\infty),E)$ denote the Banach space
endowed the maximum norm $\|u\|_{C}=\sup_{t\in[-\tau,\infty)}\|u(t)\|$. For $u\in C([-\tau,\infty),E)$ and $t\in[0,\infty)$, we denote $u_{t}\in C([-\tau,0],E)$, $u_{t}(s)=u(t+s)$, $s\in[-\tau,0]$.
 Let $\varphi\in  C([-\tau,0],E)$,
 we study the following initial value problem of the evolution equation with delay
 $$\left\{\begin{array}{ll}
u'(t)+Au(t)=F(t,u(t),u(t-\tau)),\ \ \  t\in J,\\[8pt]
u_{0}=\varphi,
 \end{array} \right.\eqno(3.17)$$
 where $-A$ generates  positive
  $C_{0}$-semigroup $T(t)(t\geq0)$ in $E$ and
  $F:J\times K\times K\to K$ be continuous.

 If there exists $u\in C([-\tau,\infty),E)$ satisfying
   $u(t)=\varphi (t)$ for $-\tau\leq t\leq0$ and
  $$u(t)=T(t)u(0)+\int^{t}_{0}T(t-s)F(s,u(s),u(s-\tau)),\ t\geq0, \eqno(3.18)$$
 then $u$ is
called a mild solution of the nonlinear initial value problem (3.17).
 Furthermore, when $\varphi\in  C([-\tau,0],K)$,
 it follows that $u(t)\geq\theta (t\in[-\tau,\infty))$
 by the characteristic of positive semigroups.

For the nonlinear initial value problem (3.17), we have the following  result.
\vskip3mm
\noindent\textbf{Lemma 3.1.} \emph{Let $E$ be an ordered Banach space
whose positive cone $K$ is normal cone,  $-A$ generate a positive
 compact semigroup $T(t)(t\geq0)$ in $E$.
Assume that $F:\R\times K\times K\rightarrow K$ is
continuous, $\varphi\in  C([-\tau,0],K)$.
If $F$ satisfies the condition (H3),
then the initial value problem (3.17) has a
 unique positive mild solution $u\in C([-\tau,\infty),K)$.}
\vskip3mm
\noindent\textbf{Proof} By the condition (H3),
 we have
$$F(t,x,y)\leq C_{1}x+C_{2}y+F(t,\theta,\theta),\ \ t\in [0,\infty),x,y\in K,\eqno(3.19)$$
namely, $\|F(t,x,y)\|\leq C_{1}\|x\|+C_{2}\|y\|+K$,
where $K=\max\limits_{t\in[0,\omega]}\|F(t,\theta,\theta)\|$.
Thus, by a standard argument as in \cite[Theorem 3.1]{Li2011},
we can prove that initial value problem (3.17) exists  positive mild solution.

Next, we show the uniqueness. Let  $u_{1},u_{2}\in C([-r,\infty),K)$
be the positive  solutions of the initial value problem (3.17),
hence they satisfy the initial value condition
 $u_{1}(t)=u_{2}(t)=\varphi(t) (-\tau\leq t\leq0)$ and (3.18).
 Let us assume that $u_{1}\leq u_{2}$, by the condition (H3),
  for every $t\geq0$, we have
\begin{eqnarray*}
&&u_{2}(t)-u_{1}(t)\\[8pt]
&=&\int_{0}^{t}T(t-s)\Big(F(s,u_{2}(s),u_{2}(s-\tau))-F(s,u_{1}(s),u_{1}(s-\tau))\Big)ds\\
&\leq&\int_{0}^{t}T(t-s)\Big(C_{1}(u_{2}(s)-u_{1}(s))+C_{2}(u_{2}(s-\tau)-u_{1}(s-\tau))\Big)ds.
\end{eqnarray*}
Define an operator $\mathcal{B}:C([-r,\infty),K)\to C([-r,\infty),K)$ as following:
$$\mathcal{B}u(t)=C_{1}u(t)+C_{2}u(t-\tau),\ \ \ t\geq0,\ u\in C([-\tau,\infty),K).$$
Clearly, $\mathcal{B}$ is a linear bounded operator with $\|\mathcal{B}\|\leq C_{1}+C_{2}$.
 Therefore,
$$u_{2}(t)-u_{1}(t)\leq\int_{0}^{t}T(t-s)\mathcal{B}(u_{2}(s)-u_{1}(s))ds,$$
which implies that
$$\|u_{2}(t)-u_{1}(t)\|\leq\int_{0}^{t}\|T(t-s)\|\cdot(C_{1}+C_{2})\cdot\|u_{2}(s)-u_{1}(s)\|ds.$$
By the Gronwall-Bellman inequality,
we have $\|u_{2}(t)-u_{1}(t)\|\equiv 0(t\geq0)$.
Hence, $u_{1}\equiv u_{2}$. \ \ $\square$
\vskip3mm
The proof of Theorem 1.3 needs the following
integral inequality of Bellman type with delay.
\vskip3mm
\noindent\textbf{Lemma 3.2.}(\cite{Li2011})\emph{Let
us assume that $\phi\in C([-r,\infty),J)$ and
there exist positive constants
$c_1, c_2,$ such that $\phi$ satisfy the integral inequality
$$\phi(t)\leq\phi(0)+c_{1}\int_{0}^{2}\phi(s)ds+c_{2}\int_{0}^{t}\phi(s-\tau)ds,\ \ t\geq0.\eqno(3.20)$$
 Then $\phi(t)\leq\|\phi\|_{C[-\tau,0]}e^{(c_{1}+c_{2}) t}$ for every $t\geq0$,
 where $\|\phi\|_{[-\tau,0]}=\max\limits_{t\in[-\tau,0]}|\phi(t)|$.}

\vskip3mm
\noindent\textbf{Proof of Theorem 1.3.}  By Theorem 1.2, the  delay
evolution equation (1.4) has a unique positive $\omega$-periodic mild solution $u^{*}\in C_{\omega}(\R,K)$. For any $\varphi\in C([-\tau,\infty),K)$,
the initial value problem (3.17)  has a unique global positive mild solution $u=u(t,\varphi)\in C([-r,\infty),K)$ by Lemma 3.1.

By the semigroup representation of the solutions, $u^{*}$ and $u$
satisfy the integral equation (3.18). Thus, by (3.18) and
assumption (H3), for any $t\geq0$, we have
 \begin{eqnarray*}\qquad\qquad u(t)-u^{*}(t)&\leq& T(t)(u(0)-u^{*}(0))
 +\int^{t}_{0}T(t-s)(C_{1}(u(s)-u^{*}(s))\\[8pt]
 &&+C_{2}(u(s-\tau)-u^{*}(s-\tau)))ds.\qquad\qquad\qquad\ \qquad(3.21) \end{eqnarray*}
 Since $T(t)(t\geq0)$ is an exponentially stable positive
$C_{0}$-semigroup, that is the growth exponent $\nu_{0}<0$, hence,
by the property of semigroup,
there is a number $M\geq1$ such that
 $$\|T(t)\|\leq Me^{\nu_{0}t},\ \ \ t\geq0.$$
 We choose the equivalent norm $|\cdot|_{0}$ by
$$|x|_{0}=\sup\limits_{t\geq0}\|e^{\nu_{0} t}T(t)x\|,$$
then $\|x\|\leq |x|_{0}\leq M\|x\|$.
Thus, we  denote the norm of $T(t)(t\geq0)$
in $(E,|\cdot|_{0}$) by $|T(t)|_{0}$ and
$|T(t)|_{0}<e^{-\nu_{0} t}$ for $t\geq0$.

 Now, by (3.21) and the normality of cone $K$ in $E$,
 we have
 \begin{eqnarray*}&&|u(t)-u^{*}(t)|_{0}\\[8pt]
 &\leq& |T(t)|_{0}\cdot|u(0)-u^{*}(0)|_{0}\\[8pt]
 &&+
 \int^{t}_{0}|T(t-s)|_{0}\Big(C_{1}|u(s)-u^{*}(s)|_{0}+C_{2}|u(s-\tau)-u^{*}(s-\tau)|_{0}\Big)ds\\[8pt]
 &\leq& e^{\nu_{0} t}|u(0)-u^{*}(0)|_{0}\\[8pt]
 &&+ \int^{t}_{0}e^{\nu_{0} (t-s)}\Big(C_{1}|u(s)-u^{*}(s)|_{0}+C_{2}|u(s-\tau)-u^{*}(s-\tau)|_{0}\Big)ds\\[8pt]
 &\leq& e^{\nu_{0} t}|u(0)-u^{*}(0)|_{0}+C_{1}e^{\nu_{0}t}\int^{t}_{0}e^{-\nu_{0}s}
 (|u(s)-u^{*}(s)|_{0})ds\\[8pt]
 &&+C_{2}e^{\nu_{0}(t-\tau)}\int^{t}_{0}e^{-\nu_{0}(s-\tau)}
 (|u(s-\tau)-u^{*}(s-\tau)|_{0})ds.\end{eqnarray*}

For \ $t\in[-\tau,\infty)$, setting\ $\phi(t)=e^{-\nu_{0}t}|u(t)-u^{*}(t)|_{0}$,
from the inequality above, it
follows that
$$\phi(t)\leq \phi(0)+C_{1}\int_{0}^{t}\phi(s)ds+C_{2}e^{-\nu_{0}\tau}\int_{0}^{t}\phi(s-\tau)ds.\eqno(3.22)$$
Hence, by Lemma 3.2, we have
$$e^{-\nu_{0}t}|u(t)-u^{*}(t)|_{0}=\phi(t)\leq C(\varphi)e^{(C_{1}+Ce^{-\nu_{0}\tau} )t},\ \ \ t\geq0,\eqno(3.23)$$
where \ $C(\varphi)=\max_{s\in[-\tau,0]}\{e^{-\nu_{0}s}|\varphi(s)-u^{*}(s)|_{0}\}$.
By the assumption (H3), $\sigma:=-\nu_{0}-(C_{1}+Ce^{-\nu_{0}\tau} )>0$,
and form (3.23) it follows that
$$|u(t)-u^{*}(t)|_{0}\leq C(\varphi)e^{-\sigma t}\to 0\ \ \ (t\to \infty).$$
Thus, the positive $\omega$-periodic solution $u^{*}$ is
globally asymptotically stable and it exponentially
attracts every positive solution of the initial value problem.
 This completes the proof of Theorem 1.3.
\hfill$\Box$

\section{Application}
In this section, we present one example, which indicates
how our abstract results can be
applied to concrete problems.
Let $\overline{\Omega}\in\R^{n}$ be a bounded domain with
a sufficiently smooth boundary $\partial\Omega$.
Let
$$A(x,D)u=-\sum^{N}_{i,j=1}a_{ij}(x)D_{i}D_{j}u+\sum_{j=1}^{N}a_{j}(x)D_{j}u+a_{0}(x)u,\eqno(4.1)$$
 be a uniformly elliptic differential operator
 in $\overline{\Omega}$, whose coefficients $a_{ij}(x), a_{j}(x)$ $(i,j=1,\cdots,n)$
 and $a_{0}(x)$ are H\"{o}der-continuous on $\overline{\Omega}$,
 and $a_{0}(x)\geq0$.
 We let $B=B(x,D)$ be a
boundary operator on $\partial \Omega$ of the form:
$$Bu:=b_{0}(x)u+\delta\frac{\partial u}{\partial\beta},\eqno(4.2)$$
where either $\delta=0$ and $b_{0}(x)\equiv1$ (Dirichlet boundary operator),
or $\delta=1$ and $b_{0}(x)\geq0$ (regular oblique derivative boundary operator;
at this point, we further assume that $a_{0}(x)\not \equiv0$ or $b_{0}(x)\not \equiv0$),
 $\beta$ is an outward pointing,
 nowhere tangent vector field on $\partial\Omega$.
Let $\lambda_{1}$ be the first eigenvalue of elliptic
operator  $A(x,D)$  under the boundary condition
 $Bu=0$. It is well known (\cite[ Theorem 1.16]{Amann76},) that $\lambda_{1}>0$.

Under the above assumptions, we discuss the existence,
uniqueness and asymptotic stability  of positive time
$\omega$-periodic solutions of the semilinear parabolic
boundary value problem
$$\left\{\begin{array}{ll}
\frac{\partial }{\partial t}u(x,t)+A(x,D) u(x,t)=f(x,t,u(x,t),u(x,t-\tau)),\ x\in \Omega, \ t\in\R,\\[10pt]
Bu=0,\quad x\in \partial\Omega,
 \end{array} \right.\eqno (4.3)$$
where\ $f:\overline{\Omega}\times\R\times\R^{2}\rightarrow\R$ a
local H\"{o}lder-continuous function which is $\omega$-periodic
in $t$, $\tau>0$ denotes the time delay.

Let\ $E=L^{p}(\Omega)(p>1)$, $K=\{u\in E|\ u(x)\geq0 \ a.e.\  x\in\Omega\}$,
then\ $E$ is an ordered Banach
space, whose positive cone $K $ is a normal regeneration cone.
Define an operator\ $A:D(A)\subset E\rightarrow E$ by:
$$D(A)= \{u\in W^{2,p}(\Omega)|\ B(x,D)u=0,\ x\in \partial\Omega\},\quad Au=A(x,D)u.\eqno(4.4)$$

If $a_{0}(x)\geq0$, then $-A$ generates an exponentially stable
analytic semigroup $T_{p}(t)(t\geq0)$
in $E$ (see \cite{Amann78}). By the maximum principle of elliptic
operators, we know that $(\lambda I+A)$ has a positive bounded
inverse operator $(\lambda I+A)^{-1}$ for $\lambda>0$,
hence $T_{p}(t)(t\geq0)$ is a positive
semigroup (see \cite{Li1996}). From the operator $A(x,D)$
has compact resolvent in $L^{p}(\Omega)$,
we obtain $T_{p}(t)(t\geq0 )$ is also a compact semigroup (see \cite{Pazy}).
Therefore, by Corollary 1.4, we have the following result.

\vskip3mm
\noindent\textbf{Theorem 4.1.} Assume that
$f:\overline{\Omega}\times\R\times\R^{2}\rightarrow\R$
is a local H\"{o}lder-continuous function which
is $\omega$-periodic in $t$ and satisfies $f(x,t,u,v)\geq0$ for \ $(x,t,u,v)\in (\overline{\Omega}\times\R\times\R^{+}\times\R^{+})$.
If the following
condition holds:

(H4) there are  constants $C_{1},C_{2}$, $C_{1}+C_{2}\in (0,\lambda_{1})$
and a function $h\in C_{\omega}(\overline{\Omega}\times\R)$ satsfying $h(x,t)\geq0$
such that
$$f(x,t,u,v)\leq C_{1}u+C_{2}v+h(x,t),\ \ \ (x,t)\in \overline{\Omega}\times\R, u,v\geq0,$$
then the delay parabolic boundary value problem (4.3)
has at least one positive $\omega$-periodic solution
$u\in C^{2,1}(\overline{\Omega}\times\R)$.
\vskip3mm

\noindent\textbf{Proof} Let$u(t)=u(\cdot,t)$,  $F(t,u(t),u(t-\tau))=f(\cdot,t,u(\cdot,t),u(\cdot,t-\tau))$,
then the delay parabolic boundary value problem
(4.3) can be reformulated as the abstract evolution equation (1.4) in $E$.
From the assumption, it is easy to see that the  conditions
of Corollary 1.4 are satisfied. By Corollary 1.4,
the delay parabolic boundary value
problem (4.3) has a time positive $\omega$-periodic
mild solution $u\in C_{\omega}(\R,E)$. By the analyticity
of the semigroup $T_{p}(t)(t\geq0)$ and the regularization method used in \cite{Amann78} , we can see that
$u\in C^{2,1}(\overline{\Omega}\times\R)$
is a classical time $\omega$-periodic solution
of the equation (4.3). This completes
the proof of the theorem. \hfill$\Box$

From Corollary 1.5 and Theorem 4.1,
we obtain the uniqueness result.

\vskip3mm
\noindent\textbf{Theorem 4.2.}  \emph{Assume that $f:\overline{\Omega}\times\R\times\R^{2}\rightarrow\R$ is a local H\"{o}lder-continuous function which
is $\omega$-periodic in $t$ and satisfies $f(x,t,u,v)\geq0$ for $(x,t,u,v)\in (\overline{\Omega}\times\R\times\R^{+}\times\R^{+})$. If the following
condition holds:
\vskip1mm
\noindent
(H5) there are constants $C_{1},C_{2}$,\ $C_{1}+C_{2}\in(0,\lambda_{1})$
such that for $y_{i},z_{i}\in K(i=1,2), y_{1}\leq y_{2},z_{1}\leq z_{2}$,
$$f(x,t,y_{2},z_{2})-f(x,t,y_{1},z_{1})\leq C_{1}(y_{2}-y_{1})+C_{2}(z_{2}-z_{1}), \ \ t\in \R,$$
\vskip1mm
\noindent then the parabolic boundary value problem (4.3)
has a unique positive $\omega$-periodic solution
$u^{*}\in C^{2,1}(\overline{\Omega}\times\R)$.}
\vskip3mm
Let $\varphi\in C(\Omega\times[-\tau,\infty))$, define a mapping
  $t\mapsto \varphi(\cdot,t)$, then we can see $\varphi\in C([-\tau,0],X)$.
  Consider the semilinear delay parabolic initial boundary value problem
$$\left\{\begin{array}{ll}
\frac{\partial }{\partial t}u(x,t)+A(x,D) u(x,t)=f(x,t,u(x,t),u(x,t-\tau)),\ x\in \Omega, \ t\geq0,\\[8pt]
Bu=0,\quad x\in \partial\Omega,\\[8pt]
u(x,t)=\varphi(x,t), (x,t)\in\Omega\times[-\tau,0],
 \end{array} \right.\eqno (4.5)$$
From Lemma 3.1, we can obtain the following existence and uniqueness results.

\vskip3mm
\noindent{\textbf{Lemma 4.1.}} \emph{Let $f:\overline{\Omega}\times\R\times\R^{2}\rightarrow\R$ is a local H\"{o}lder-continuous function, and $f(x,t,u,v)\geq0$ for every $(x,t,u,v)\in (\overline{\Omega}\times\R\times\R^{+}\times\R^{+})$.
If the following condition
\vskip1mm
 \noindent  (H6) there are constants $C_{1},C_{2}$ satisfying  $C_{1}+C_{2}e^{\lambda_{1}\tau}\in(0,\lambda_{1})$,
 such that for any $y_{i},z_{i}\in K(i=1,2)$ with $ y_{1}\leq y_{2},z_{1}\leq z_{2}$,
 $$f(x,t,y_{2},z_{2})-f(x,t,y_{1},z_{1})\leq C_{1}(y_{2}-y_{1})+C_{2}(z_{2}-z_{1}), \ \ (x,t)\in \Omega\times\R,$$
\noindent holds, then the delayed parabolic initial boundary
value problem\ (4.5) has a uniqueness positive solution
 $u\in C([-\tau,\infty),L^{p}(\overline{\Omega}))\cap C^{2,1}(\overline{\Omega}\times(0,\infty))$.}

\vskip2mm
Hence, From Corollary 1.6, we can derive the asymptotic
stability of the positive $\omega$-periodic solution for
the delay parabolic boundary value problem (4.3).

\vskip3mm
\noindent{\bf Theorem 4.3}\emph{ Assume that
$f:\overline{\Omega}\times\R\times\R^{2}\rightarrow\R$
is a local H\"{o}lder-continuous function which
is $\omega$-periodic in $t$ and satisfies $f(x,t,u,v)\geq0$
for $(x,t,u,v)\in (\overline{\Omega}\times\R\times\R^{+}\times\R^{+})$.
If the condition (H6) holds,
then the unique positive $\omega$-periodic solution
$u^{*}\in C^{2,1}(\overline{\Omega}\times\R)$ of equation (4.3)
is globally asymptotically stable and it exponentially
attracts every positive solution of the initial value problem.
}
\\

\bibliographystyle{abbrv}

\begin{thebibliography}{99}
\small
\setlength{\parskip}{0pt}
\setlength{\baselineskip}{12pt}
\vspace{4pt}

\bibitem{Ahmed01}
N. Ahmed, Measure solutions of impulsive evolution differential inclusions and optimal control,   Nonlinear Anal., 47(2001): 13-23.


\bibitem{Amann76} H. Amann, Nonlinear operators in ordered Banach spaces and some applicitions to nonlinear boundary value problem, In: Nonlinear operators and the Calculus of Variations, Lecture Notes in Mathmematics, Springer-Verlag, Berlin and New YorK, 1976, pp. 1-55.
 \bibitem{Amann78}
H. Amann, Periodic solutions of semilinear parabolic equations, In: L. Cesari, R. Kannan, R. Weinberger (Eds.),
 Nonlinear Anal. A Collection of Papers in Honor of Erich H. Rothe, New York: Academic Press, 1978: 1-29.

\bibitem{Burton1991} T. Burton, B. Zhang, Periodic solutions of abstract differential equations with infinite delay. J. Differential Equations 90(1991):357-396.

 \bibitem{Deimling}  K. Deimling, Nonlinear Functional Analysis, Springer, New York (1985).


\bibitem{Hale1993} J. Hale, S. Lunel, Introduction to Functional-differential equations, Applied

\bibitem{Kpoumie} M. EK. Kpoumi\`{e}, K. Ezzinbi, D. B\'{e}koll\`{e}, Periodic Solutions for Some Nondensely Nonautonomous
Partial Functional Differential Equations in Fading
Memory Spaces, Differ. Equ. Dyn. Syst.,(2016):1-21.

\bibitem{Li1996}  Y. Li, The positive solutions of abstract semilinear evolution equations and their applications,Acta Math. Sin., 39 (1996) :666-672 (in Chinese).

\bibitem{Li98}Y. Li, Periodic solutions of semilinear evolution equations in Banach spaces. Acta Math. Sin. 41(1998) 629-636 (in Chinese).

\bibitem{Li2005x}Y. Li,  Existence and uniqueness of positive periodic solution for abstract semilinear evolution equations, J. Syst. Sci. Math. Sci., 25( 2005): 720-728  (in Chinese).


\bibitem{Li2009} Y. Li, Existence and uniqueness of periodic solution for a class of semilinear evolution equations, J. Math. Anal. Appl. 349(2009): 226-234.

\bibitem{Li2011}
Y. Li, Existence and asymptotic stability of
 periodic solution for evolution equations with delays,
 J. Funct. Anal. 261 (2011):1309-1324.

\bibitem{QLi2015}  Q. Li, Y. Li,  Existence of positive periodic solutions for
abstract evolution equations,  Advances in Difference Equations, (2015) 2015:135, 12 pages.


\bibitem{Liang2017A}J. Liang et al, Periodicity of solutions to the Cauchy problem
for nonautonomous impulsive delay evolution
equations in Banach spaces,  Anal. Appl., 15(2017):457-476.

\bibitem{Liang2017D} J. Liang, J.H. Liu,  T.J. Xiao, Condensing operators and periodic solutions of infinite delay impulsive evolution equations, Discrete and Continuous Dynamical Systems - Series,
 10(2017):475-485.

\vskip2mm
\bibitem{Liu1998} J. Liu, Bounded and periodic solutions of finite delays evolution equations, Nonlinear Anal. 34 (1998) 101-111.
\vskip2mm
\bibitem{Liu2000} J. Liu, Periodic solutions of infinite delay evolution equations,  J. Math. Anal. Appl. 247(2000): 644-727.
\vskip2mm
\bibitem{Liu2003} J. Liu, Bounded and periodic solutions of infinite delay evolution equations . J. Math. Anal. Appl. 2003, 286(2003):705-712.

 \bibitem{Pazy} A. Pazy. Semigroups of Linear Operators and Applications to Partial
Differential Equations. Berlin: Springer-Verlag, 1983.




\bibitem{ShiLiSun11} H. Shi,  W. Li, H. Sun, Existence of mild solutions for abstract mixed
type semilinear evolution equation , Turk. J. Math. 35(2011): 457-472.

\bibitem{Triggiani} R. Triggiani, On the stabilizability problem in Banach space, J. Math. Anal. Appl. 52(1975):383-403.


\bibitem{Wang2016} R. Wang, P. Zhu, New Results on Periodic Solutions to Impulsive Nonautonomous Evolutionary Equations with Time Delays, J. Math. Sci., (2016) 212: 412.

\bibitem {Wu1996}J. Wu, Theory and Applications of Partial Functional Differential Equations, Appl. Math. Sciences,  vol. 119, New York: Springer, 1996.

\bibitem{Xiang1992}X. Xiang, N.U. Ahmed, Existence of periodic solutions of semilinear evolution equations with time lags, Nonlinear
Anal., 18 (1992) :1063-1070.

\bibitem{Yang} H. Yang, Q. Li, Asymptotic stability of positive periodic
solution for semilinear evolution equations, Advances in Difference Equations 2014(2014):197,10 pages.



\end{thebibliography}

\end{document}